\documentclass[12pt]{article}
\usepackage{latexsym}
\usepackage{amsfonts}
\usepackage{epsfig}
\newtheorem{theorem}{Theorem}

\def\C{{\mathbf{C}}}
\def\bC{{\mathbf{\overline{\mathbf{C}}}}}
\def\R{{\mathbf{R}}}

\begin{document}
\title{Spectral loci of Sturm--Liouville operators with polynomial potentials}
\author{Alexandre Eremenko
and Andrei Gabrielov\thanks{Both authors are supported by
NSF grant DMS-1067886.}}
\maketitle
\begin{abstract}
We consider differential equations $-y^{\prime\prime}+P(z,a)y=\lambda y$,
where $P$ is a polynomial of the independent variable $z$ depending on
a parameter $a$. The spectral locus is the set of $(a,\lambda)$ such that
the equation has a non-trivial solution tending to zero on two fixed rays
in the complex plane. We study the topology of the spectral loci
for polynomials $P$ of degree $3$ or $4$
with respect to $z$.

MSC: 81Q05, 34M60, 34A05.

Keywords: one-dimensional Schr\"odinger operators,
quasi-exact solvability, PT-symmetry.
\end{abstract}

We consider differential operators
$$L: y\mapsto -y^{\prime\prime}+Py,$$
where the potential $P$ is a polynomial of degree $d$.
When $d\in \{0,1,2\}$, the general solution of $Ly=0$ can be expressed
in terms of special functions
(elementary, Airy or Weber functions, respectively). The eigenvalue problem
for $d=2$ (harmonic oscillator) plays an important role in quantum mechanics.

We mostly consider the cases $d=3$ and $d=4$.

Cubic and quartic oscillators were studied extensively from the very
beginning of quantum mechanics, mostly by perturbative methods.
Cubic oscillator arises in quantum field
theory \cite{Z} and in the theory of Painlev\'e equations \cite{M1,M2}.

We normalize the equation by an affine transformation
of the independent variable to obtain:
$P(z)=z^d+O(z^{d-2}), \; z\to\infty.$

Every solution $y$ of $Ly=\lambda y,$ for $\lambda\in\C$
is an entire function.
Consider the Stokes sectors 
$$S_j=\left\{ z: \left|\arg z-\frac{2\pi j}{d+2}\right|<\frac{\pi}{d+2}\right\},
 \; 0\leq j\leq d+1.$$
For each $j\in\{0,\ldots,d+1\}$
a solution $y$ either
grows exponentially along all rays from the origin in $S_j$
or tends to zero exponentially
along every such ray in $S_j$.

We choose two non-adjacent sectors,
and impose the boundary condition that a solution $y$ of
$Ly=\lambda y$  tends to $0$ in
these sectors. Those $\lambda$ for which such $y\neq 0$ exists
are called eigenvalues, and the corresponding $y$'s eigenfunctions. 

The interest of physicists to such eigenvalue problems originates from the
paper \cite{BW} where a systematic study of analytic continuation
of eigenvalues as parameters vary in the complex plane was made for
the first time.

With such boundary conditions, the problem has an infinite discrete
spectrum and eigenvalues tend to infinity \cite{Sib}. To each eigenvalue
corresponds a one-dimensional eigenspace.

Suppose that the polynomial potential depends analytically
on a parameter $a\in\C^n$. Then the spectral locus $Z$ is defined as the set
of all pairs $(a,\lambda)\in\C^{n+1}$ such that the differential equation
$$-y^{\prime\prime}+P(z,a)y=\lambda y$$
has a solution $y$ satisfying the boundary conditions.
The spectral locus is an analytic hypersurface in $\C^{n+1}$; it is the zero-set
of an entire function $F(a,\lambda)$ which is called the spectral determinant
\cite{Sib}.

The multi-valued function $\lambda(a)$ defined by
$F(a,\lambda)=0$
has the following property: its only singularities are algebraic ramification
points, and there are finitely many of them over every compact set in
the $a$-space \cite{EG1}.

Next we discuss connectedness of the spectral locus.

\begin{theorem} $\cite{EGx}$ For the cubic oscillator

\begin{equation}\label{cub}
-y^{\prime\prime}+(z^3-az)y=-\lambda y,\quad y(\pm i\infty)=0, 
\end{equation}
the spectral locus is a smooth irreducible curve in $\C^2$.
\end{theorem}
\begin{theorem} $\cite{EG1}$ For the even quartic oscillator
$$-y^{\prime\prime}+(z^4+az^2)y=\lambda y,\quad y(\pm\infty)=0,$$
the spectral locus consists of two disjoint smooth irreducible curves
in $\C^2$, one corresponding to even eigenfunctions, another to odd ones.
\end{theorem}

These theorems can be generalized to polynomials of arbitrary degree
if we use all coefficients as parameters \cite{H,A,AG}.

However if we consider a subfamily of all polynomials
of given degree, then
the spectral
locus can be reducible in an interesting way \cite{BB}.

One example is the family of quasi-exactly solvable quartics $L_J$
\begin{equation}\label{quart}
-y^{\prime\prime}+(z^4-2bz^2+2Jz)y=\lambda y,
\quad y(re^{\pm\pi i/3})\to0,\quad r\to+\infty.
\end{equation}
When $J$ is a positive integer, this problem has $J$
elementary eigenfunctions
of the form
$p(z)\exp(z^3/3-bz)$, with a polynomial $p$.
The $(b,\lambda)$ corresponding to these
eigenfunctions form the quasi-exactly solvable part $Z_J^{QES}$
of the spectral locus $Z_J$, and $Z_J^{QES}$ is an algebraic curve.
\begin{theorem} $\cite{EGx}$ For every positive integer $J$,
the QES spectral locus $Z_J^{QES}$ is a smooth irreducible curve in $\C^2$.
\end{theorem}
Thus the curve $Z_J^{QES}$ is an irreducible component of $Z_J$.
It is not known whether $Z_J\backslash Z_J^{QES}$ is irreducible.

A similar phenomenon occurs in degree $6$: there are one-parametric families
of even quasi-exactly solvable sextics \cite{Tu,U}, and for each
such family the
quasi-exactly solvable part of the spectral locus is a smooth irreducible
algebraic curve \cite{EGx}.

When $J\to\infty$, an appropriate rescaling of $Z_J^{QES}$ tends
to the spectral locus of one-parametric cubic family (\ref{cub}), and
a rescaling of the sextic
QES spectral locus tends to the spectral locus of the even quartic family
\cite{U}.

Now we consider Hermitian and PT-symmetric operators.
An eigenvalue problem may be
invariant under a reflection with respect to a line in the complex $z$-plane.
Without loss of generality, we can take this line to be the real line, and
the reflection to be the complex conjugation. Then the coefficients
of the potential are real. For the boundary conditions, two cases are possible:

a) Each of the two boundary conditions is preserved by the reflection.
In this case the
problem is Hermitian.

b) The two boundary conditions are interchanged by the reflection. Such problems
are called PT-symmetric. (Physicists prefer to choose the 
reflection with respect to the imaginary line in this case. PT stands
for ``parity and time''.)

For example, there is a one-parametric family of PT-symmetric cubics
(\ref{cub}),
with real parameter $a$. 
There are two different two-parametric families of
PT-symmetric quartics which we call type I and II:
\begin{equation}\label{1}
-y^{\prime\prime}+(-z^4+az^2+cz)y=-\lambda y, \quad y(\pm i\infty)=0,
\end{equation}
and
\begin{equation}
\label{2}
-y^{\prime\prime}+(z^4-2bz^2+2Jz)y=\lambda y,\quad y(re^{\pm\pi i/3})\to 0,
\quad r\to+\infty.
\end{equation}
The real spectral locus $Z(\R)$ is the subset of the spectral locus $Z$
which consists of points with real coordinates.

We begin with the cubic PT-symmetric spectral locus (\ref{cub}).

\begin{theorem}\label{t3} $\cite{EG2}$ For every integer $n\geq 0$,
there exists a simple curve $\Gamma_n\subset\R^2$,
that is the image of a proper analytic embedding of a line,
and which has these properties:
\newline
(i) For every $(a,\lambda)\in\Gamma_n$ problem (\ref{cub})
has an eigenfunction with $2n$ non-real zeros.
\newline
(ii) The curves $\Gamma_n$ are disjoint and
the real spectral locus of (\ref{cub})
is $\bigcup_{n\geq 0}\Gamma_n$
\newline
(iii) The map
$$\Gamma_n\cap\{(a,\lambda):a\geq 0\}
\to \R_{\geq 0},\quad
(a,\lambda)\mapsto a$$ is a $2$-to-$1$ covering.
\newline
(iv) For $a\geq 0$, $(a,\lambda)\in\Gamma_n$ and
$(a,\mu)\in\Gamma_{n+1}$ imply $\mu>
\lambda.$
\end{theorem}
The following computer-generated plot of the real
spectral locus of (\ref{cub}) is taken from Trinh's thesis
\cite{T}. Theorem \ref{t3} rigorously establishes
some features of this picture. At the points whose abscissas are marked
in Fig. 1, pairs of real eigenvalues collide and escape to the complex plane.
Theorem 4 proves rigorously that infinitely many such points exist.
\bigskip
\begin{center}
\epsfxsize=4.5in%
\centerline{\epsffile{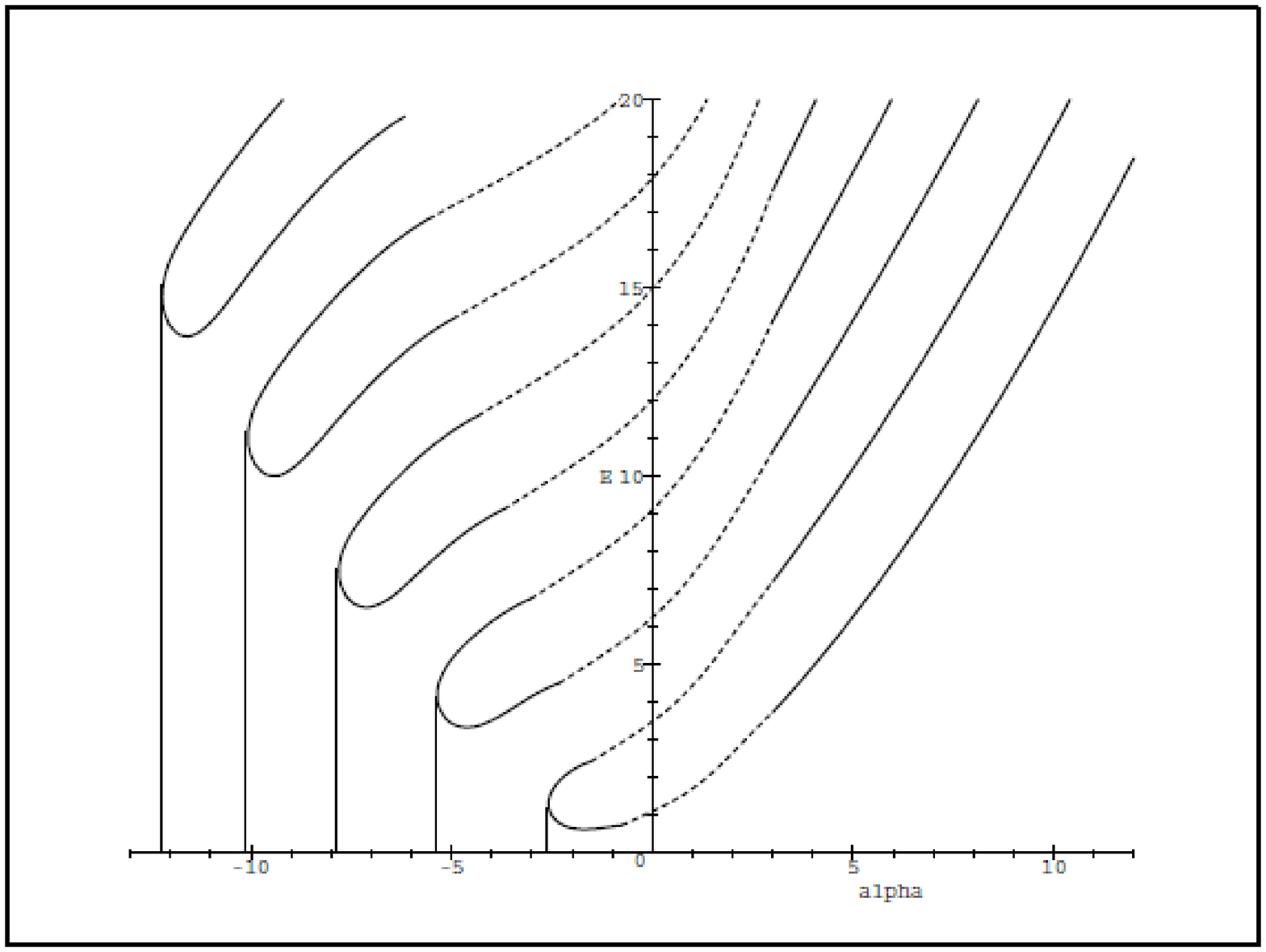}}
\medskip
\nopagebreak
\noindent Fig.~1. Real spectral locus for
$PT$-symmetric cubic.
\end{center}
Consider the PT-symmetric quartic family (\ref{1}) of type I:
It is equivalent to the $PT$-symmetric family
$$-w^{\prime\prime}+(z^4+az^2+icz)w=\lambda w,\quad w(\pm\infty)=0,$$
studied by Bender, {\em et al} \cite{BB} and Delabaere and Pham
\cite{DP}.

\begin{theorem}\label{t4} $\cite{EG2}$ The real spectral locus of
(\ref{1}) consists of disjoint smooth
analytic properly
embedded surfaces $S_n\subset\R^3, \; n\geq 0,$
homeomorphic to a punctured disk. For $(a,c,\lambda)\in S_n$,
the eigenfunction has exactly $2n$ non-real zeros.
For large $a$, projection of $S_n$ on the $(a,c)$ plane
approximates the region $9c^2-4a^3\leq 0$.
\end{theorem}
Numerical computation suggests that the surfaces have the
shape of infinite funnels with the sharp end stretching
towards $a=-\infty,\,c=0$, and that the section of $S_n$
by every plane $a=a_0$ is a closed curve.

Theorem \ref{t4} implies that this section is compact for
large $a_0$.
The following computer-generated plot is taken from Trinh's thesis:
\begin{center}
\epsfxsize=4.5in%
\centerline{\epsffile{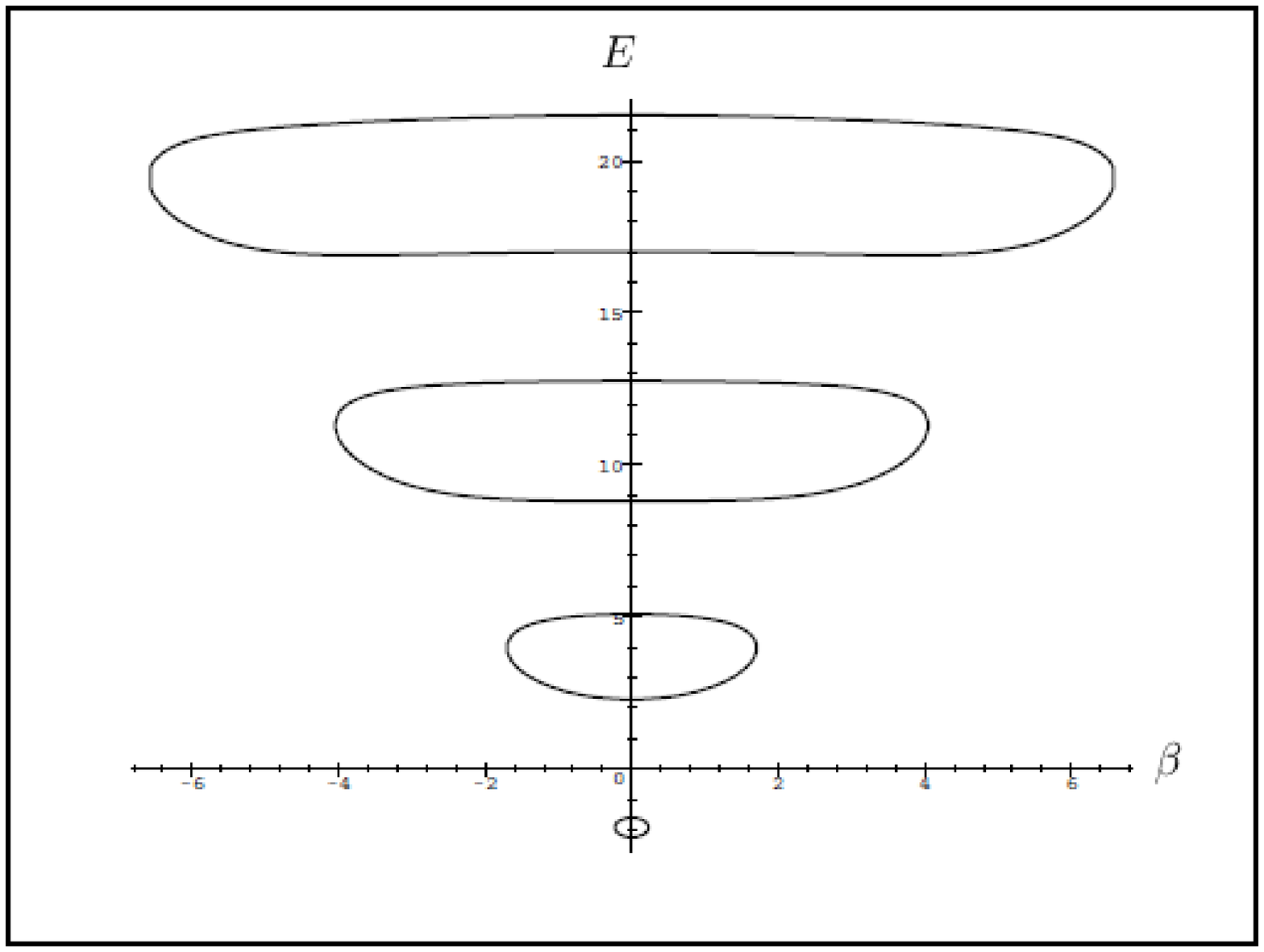}}
\nopagebreak
\noindent Fig.~2. Section of the surfaces $S_0,\ldots,S_3$
by the plane $a=-9$.
\end{center}

The PT-symmetric quartic family of the type II is more complicated,
due to the presence of the QES spectrum.
Let $Z_{J}^{QES}(\R)$ be the real QES spectral locus of the operator
$L_{J}$, defined in (\ref{quart}).
\begin{theorem} $\cite{EG4}$
For $J=n+1>0$, $Z_{n+1}^{QES}(\R)$ consists of $[n/2]+1$ disjoint
analytic curves $\Gamma_{n,m},\;0\leq m\leq[n/2].$

For $(b,\lambda)\in\Gamma_{n,m}$, the eigenfunction has $n$ zeros, $n-2m$
of them real.

If $n$ is odd, then $b\to+\infty$ on both ends of $\Gamma_{m,n}$.
If $n$ is even, the same holds for $m<n/2$, but on the ends of
$\Gamma_{n,n/2}$ we have $b\to\pm\infty$.

If $(b,\lambda)\in\Gamma_{n,m}$ and $(b,\mu)\in\Gamma_{n,m+1}$ and
$b$ is sufficiently large, then $\mu>\lambda$.
\end{theorem}  
It follows from these theorems that in each family, there are infinitely
many parameter values where pairs of real eigenvalues collide and escape 
from the real line to the complex plane.

In the quartic family of the type II, another interesting feature
of the real spectral locus is present: for some parameter values
the QES spectral locus crosses the rest of the spectral locus.
This is called ``level crossing''.

\begin{theorem} $\cite{EG3}$ The points $(b,\lambda)\in Z_{J}^{QES}$ where the level
crossing occurs are the intersection points of $Z_{J}^{QES}$ with
$Z_{-J}$. For each $J\geq 1$ there are infinitely many such points,
in general, complex. When $J$ is odd, there are infinitely many
level crossing points with $b_k<0$ and real $\lambda_k$. We have 
$$b_k\sim-((3/4)\pi k)^{2/3},\; k\to\infty.$$
\end{theorem}
The only known general result of reality of eigenvalues of $PT$-symmetric
operators is a theorem of K. Shin \cite{Shin}, which for our quartic of type II
implies that all eigenvalues are real if $J\leq 0$.

We have the following extensions of this result.
\begin{theorem} $\cite{EG5}$
For every positive integer $J$, all non-QES eigenvalues of $L_{J}$
are real.
\end{theorem}
\begin{theorem} $\cite{EG5}$
All eigenvalues of $L_{J}$ are real for every real $J\leq 1$
(not necessarily integer).
\end{theorem}

\begin{center} \epsfxsize=5.0in%
\centerline{\epsffile{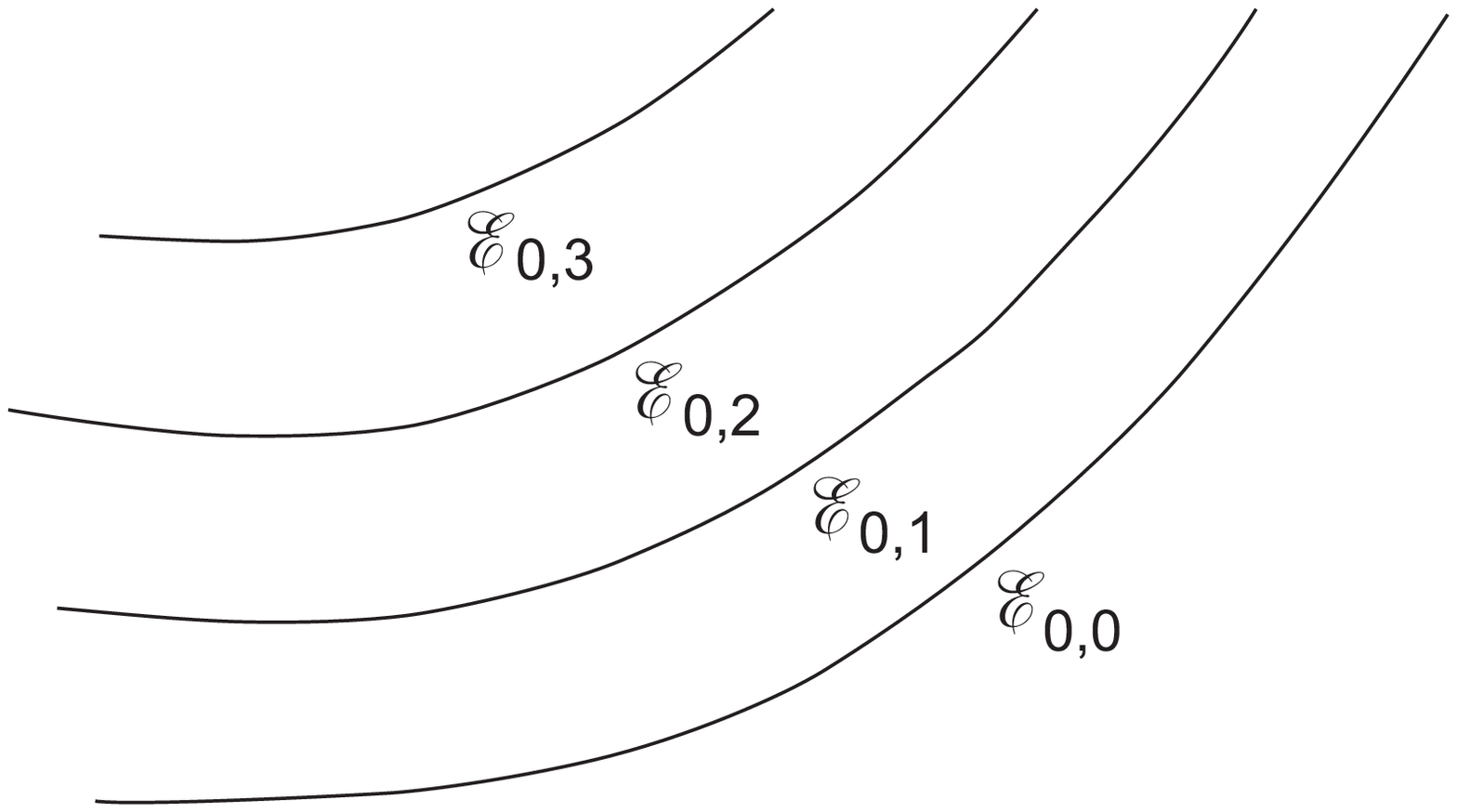}}
\nopagebreak\medskip
\noindent Fig. 3. $Z_0(\R)$.

\epsfxsize=5.0in%
\centerline{\epsffile{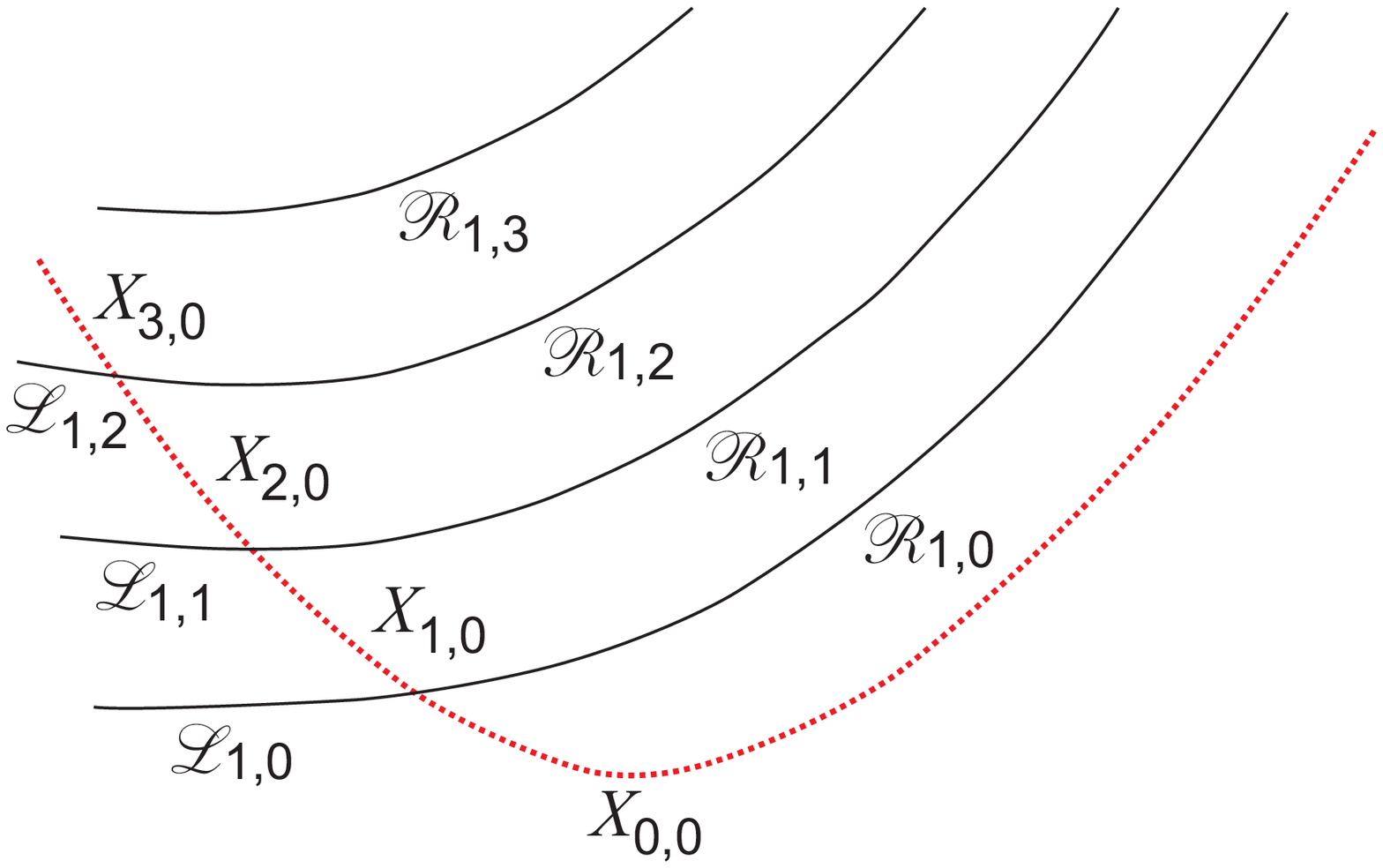}}
\nopagebreak\medskip
\noindent Fig. 4. $Z_{1}(\R).$
\end{center}

\begin{center}
\epsfxsize=5.0in%
\centerline{\epsffile{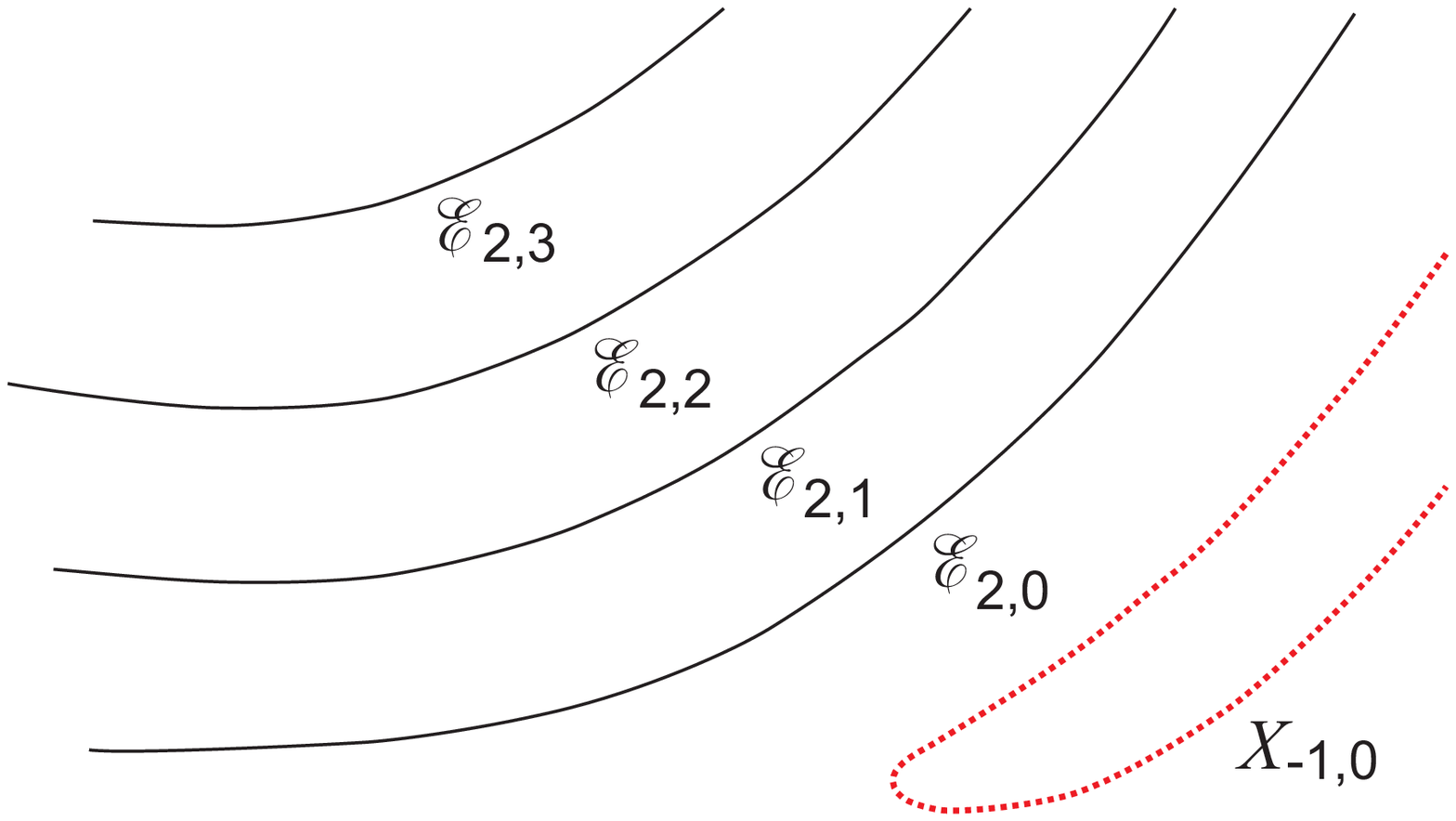}}
\nopagebreak\medskip
\noindent Fig. 5. $Z_2(\R)$.

\medskip

\epsfxsize=5.0in%
\centerline{\epsffile{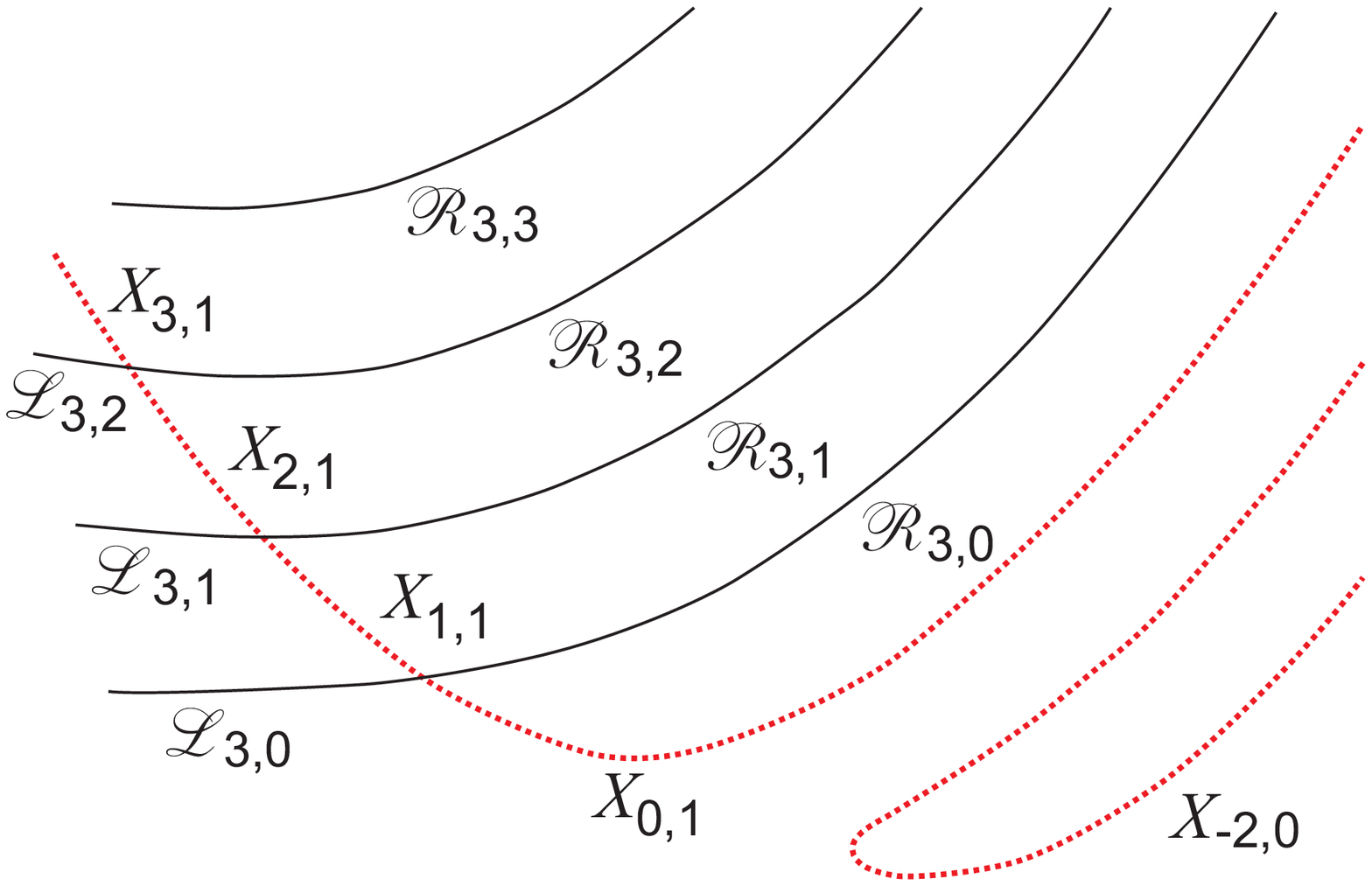}}
\nopagebreak\medskip
\noindent Fig. 6. $Z_3(\R).$
\end{center}

We briefly describe the methods of proofs. There are three main ingredients:
\vspace{.1in}

\noindent
a) Nevanlinna parametrization of the spectral locus.
\vspace{.1in}

\noindent
b) Asymptotics at infinity (singular perturbation theory).
\vspace{.1in}

\noindent
c) Darboux transform of the QES quartic.
\vspace{.1in}

Nevanlinna parametrization is the following.
Let $Z$ be the spectral locus of the problem
$$-y^{\prime\prime}+P(z,a)y=\lambda y,\quad y(z)\to 0,\quad z\in S_j\cup S_k,$$
where $S_j$ and $S_j$ are non-adjacent Stokes sectors.
Let $(a,\lambda)\in Z$, and $y_0$ an eigenfunction.
Let $y_1$ be a second linearly independent solution. Then $f=y_0/y_1$
satisfies the Schwarz differential equation
$$\frac{f^{\prime\prime\prime}}{f^\prime}-
\frac{3}{2}\left(\frac{f^{\prime\prime}}{f^\prime}\right)^2=
-2(P-\lambda).$$
This function $f$ is meromorphic in $\C$, has no critical points
and has $d+2$ asymptotic values (radial limits along rays),
one asymptotic value in each Stokes sector.
Asymptotic values in $S_j$ and $S_k$ are $0$. Asymptotic values
in adjacent sectors are distinct.

In the opposite direction: if we have a meromorphic function in $\C$
without critical points and with finitely many asymptotic values, then
it satisfies a Schwarz equation whose RHS is a polynomial \cite{N}.
The degree of this polynomial is the number of asymptotic tracts minus $2$.
(Asymptotic tracts are in one-to-one correspondence with the logarithmic
singularities of the inverse function).

Asymptotic values are meromorphic functions on $Z$ which serve as local
parameters. These are called the Nevanlinna parameters;
they are simply related
to the Stokes multipliers of the linear ODE \cite{Sib,M1,M2}.

Functions $f$ of the above type with given set of asymptotic values
$A=\{ a_0,\ldots,a_{d+1}\}$ have the property that
$$f:\C\backslash f^{-1}(A)\to\bC\backslash A$$
is a covering map. For a fixed $A$ such covering map can be completely 
described by certain combinatorial information, a cell decomposition
of the plane. These cell decompositions label the charts of our
description of the spectral locus.

It is important that we know exactly which cell decompositions can occur
and how the cell decomposition changes when the point $(a_0\ldots,a_{d+1})$
describes a closed loop in $\C^{d+2}$.
This induces an action of the braid group on the set of special
cell decompositions of the plane which can be explicitly computed.

This reduces the problem of parametrization of a spectral
locus to combinatorial topology.

For QES operators we use the Darboux transform.
Let $-D^2+V$ be a second order linear differential operator with
potential $V$. Let $\phi_0,\ldots\phi_n$ be some eigenfunctions
with eigenvalues $\lambda_0,\ldots,\lambda_n$. The
transformed operator is
$$-D^2+V-2\frac{d^2}{dz^2}\log W(\phi_0,\ldots,\phi_n),$$
where $W$ is the Wronski determinant. The eigenvalues of the
transformed operator are exactly those eigenvalues of $-D^2+V$ which
are {\em distinct} from $\lambda_0,\ldots,\lambda_n$.

We use the Darboux transform to kill the QES part
of the spectrum of $L_J$ and it turns out that the transformed
operator is $L_{-J}$!
Our study of the QES locus of the quartic family gives the following
interesting identities with elementary intergals.

To describe these identities, we address the question studied
by Heine and Stieltjes in XIX century.
Let $h$ and $p$ be polynomials. When does $y=pe^h$ satisfy
a linear differential equation $y^{\prime\prime}+Py=0$
with a polynomial $P$? 

\begin{theorem} The following conditions are equivalent:

a) $p^{\prime\prime}+2p^\prime h^\prime$ is divisible by $p$,

b) $p^{-2}e^h$ has no residues,

c) zeros of $p$ satisfy the system of equations
$$\sum_{j:j\neq k}\frac{1}{z_k-z_j}=-h^\prime(z_k),\; 1\leq k\leq\deg p.$$
\end{theorem}
Now take $h(z)=z^3/3-bz.$

\begin{theorem} Let $p$ be a polynomial. All residues of $y=p^{-2}e^{-2h}$
vanish if and only if there exists a constant $C$ and a polynomial $q$
such that
$$\left(p^2(-z)-\frac{C}{p^2(z)}\right)e^{-2h(z)}=\frac{d}{dz}\left(
\frac{q(z)}{p(z)}e^{-2h(z)}\right).$$
Moreover, if this happens then
$$C=(-1)^n2^{-2n}\frac{\partial}{\partial\lambda}Q_{n+1},$$
where $\lambda=y^{\prime\prime}/y-z^4+2bz^2-2(n+1)z$, and $Q_{n+1}(b,\lambda)=0$
is the equation of the QES spectral locus of $L_{n+1}$.
\end{theorem}

This was conjectured in \cite{EG3} on the basis of calculations
with Darboux transform
of $L_{n+1}$ and proved by E. Mukhin and V. Tarasov \cite{MT}.

Purdue University

www.math.purdue.edu/\~{}eremenko

www.math.purdue.edu/\~{}agabriel
\end{document}